\documentclass[a4paper,10pt]{article}
\usepackage[utf8x]{inputenc}
\usepackage{amssymb}

\title{Lagrangian shadows of ample algebraic divisors}
\author{Nikolay A. Tyurin\footnote{BLTPh (Dubna) and  NRU HSE (Moscow), ntyurin@theor.jinr.ru. The article was prepared within the framework of a  subsidy granted to the HSE by the Government
         of the Russian Federation for the implementation of the Global 
 Competitiveness Program}}

\begin{document}

\maketitle

\begin{abstract} In the framework of Special Bohr - Sommerfeld geometry it was established that an ample divisor in compact algebraic variety can define
almost canonically certain real submanifold which is lagrangian with respect to the corresponding Kahler form. It is natural to call it ``lagrangian shadow'';
below we emphasize this correspondence and present some simple examples, old and new. In particular we show that for irreducible divisors from the linear system
$\vert - \frac{1}{2} K_{F^3} \vert$ on the full flag variety $F^3$ their lagrangian shadows are Gelfand - Zeytlin type lagrangian 3 - spheres.

\end{abstract}

In preprint [1] one proposes a new programme which was called Special Bohr - Sommerfeld geometry, applicable in the broadest context to
compact simply connected symplectic manifolds with integer symplectic forms. The most interesting particular case is given by the consideration
of a compact simply connected algebraic variety $X$ with an ample line bundle $L$. Then there exists a Kahler form $\omega$ on
$X$  with the integer cohomology class, equals to  $c_1(L) \in H^2(X, \mathbb{Z})$ so after fixing a hermitian structure on $L$
one can apply Special Bohr - Sommerfeld geometry and get the moduli space of special Bohr - Sommerfled lagrangian cycles, see [1].
Recall that this set is formed by Bohr - Sommerfeld lagrangian cycles which are special with respect to holomorphic sections of
$L$.  

Now suppose that we have an ample divisor $D \subset X$ in a compact simply connected algebraic variety $X$ of complex dimension $n$. By the very definition,
see [2], the complete linear system $\vert kD \vert$ for certain $k \in \mathbb{N}$ gives an embedding $\phi: X \to \mathbb{C} \mathbb{P}^N$
to the projective space. Fixing a standard Fubini - Study metric on $\mathbb{C} \mathbb{P}^N$ one gets a Kahler metric of the Hodge type
on $X$ such that the Kahler form $\omega_k$ has the cohomology class Poincare dual to $k [D] \in H_{2n-2}(X, \mathbb{Z})$; and at the same time
it gives a hermitian structure $\vert . \vert_k$ on the line bundle $L_{kD} = \phi^* {\cal O}(1) \to X$ which corresponds to the standard hermitian structure  on
${\cal O}(1) \to \mathbb{C} \mathbb{P}^N$ coming after the fixing of the Kahler structure on the last projective space. 

Consider the holomorphic section $h_D \in H^0(X, L_{kD})$ with the multiple zero locus $(h_D)_0 = k D$; this section is defined up to 
scaling, but the following real function
$$
\psi_D^k = - \rm{ln} \vert h_D \vert_k
$$
is correctly defined on the complement $X \backslash D$. In the last expression for the function we indicate the level $k$ since it is possible
to consider the embeddings of different degrees and the result of our construction evidently depends on the degree.

Since the Kahler metric with the Kahler form $\omega_k$ has been fixed on $X$ one can consider the gradient flow generated by the function $\psi^k_D$
on the complement $X \backslash D$. Since $D$ is the region where $\phi^k_D$ goes to $+ \infty$ the behavior of the gradient flow near $D$ is clear, so $D$ can be understood 
as ``infinite maximum set'' for the function $\psi^k_D$. 

In [1] it was shown that a Bohr - Sommerfeld lagrangian submanifold $S \subset X \backslash D$ is special w.r.t. the holomorphic section $h_D$
iff $S$ is stable with respect to the gradient flow generated by $\psi^k_D$:
$$
F^t_{\psi_D^k} (S) = S \quad \rm{for} \quad \rm{any} \quad t.
$$

 Therefore in this situation we are interested in other ``finite'' critical points of $\psi^k_D$ and in the ``final'' trajectories of the gradient
flow $F^t_{\psi^k_D}$. At a non critical point the trajectories cann't form a local submanifold of dimension greater than $n$: according to Milnor, see [3], every 
``finite'' critical point of $\psi^k_D$  must have
the Morse index less or equal to $n$ (since $\psi^k_D$ is strictly convex on $X \backslash D$ and plurisubharmonic w.e.t. $\omega_k$, see [1]).
On the other hand, this local submanifold must be isotorpical due to the same reason, therefore if $\psi^k_D$ admits
``finite'' critical points of Morse index $n$ then the ``finite'' trajectories give us some {\it lagrangian} patches, and in good
cases one can forms lagrangian submanifolds or lagrangian cycles from the set of these patches. 

We would like to call such a submanifold or cycle {\bf lagrangian shadow of level $k$} for a given ample algebraic divisor $D$ and denote is as $Sh^{Lag}(D)$.
 Indeed, it is defined almost uniquelly:
we have fixed in the construction above only a standard Fubini - Study metric on the projective space, but the shadows in the real life as well are not unique ---
they depend on the Sun's position in the sky (in particular sometimes shadows vanish...). The main interest in the presented construction follows from the fact
that we start with pure algebraic situation (algebraic variety, ample divisor) and get something from other part of the realm of Geometry --- lagrangian.

Today this correspondence has just a conjectural explanation, therefore we are speaking about it as phenomenological
observations. The first two examples have been presented in [1]:

{\bf Example 0.} Let $X$ be $\mathbb{C} \mathbb{P}^1$ with the standard Fubuni -Study metric and the degree of $D$ is 2. Then if $D$ is an irreducible element of $\vert 2h \vert$ then
its lagrangian shadow $Sh^{Lag}(D) = S^1$ is a circle; if $D$ is reducible, then the shadow doesn't exist (see [1]).

{\bf Example 1.} Let $X$ be the complex quadric $\mathbb{C} \mathbb{P}^1 \times \mathbb{C} \mathbb{P}^1$ with the standard Kahler structure,
and $D$ has bi - degree (1,1) (so $L_D = {\cal O} (1,1)$). Then if $D$ is irreducible then $Sh^{Lag}(D) = S^2$, and this lagrangian sphere
is Hamiltonian isotopical to the antidiagonal embedding $S^2 = \{ [x_0: x_1] \times [\bar x_0: \bar x_1] \}$; if $D$ is reducible then
the lagrangian shadow doesn't exist (see [1]). In this example we have homological non trivialty of the lagrangian shadow which present
in this case the class $(h, -h)$ up to sign.

{\bf Example 1a.} Generalizng the previous example, take $X = \mathbb{C} \mathbb{P}^n \times \mathbb{C} \mathbb{P}^n$ with the standard product Kahler structure.
Let $D$ be  the zeroset of a holomorphic section of the bundle ${\cal O}(1,1)$ then again if $D$ is irreducible then $Sh^{Lag}(D)$ is isomorphic
to $\mathbb{C} \mathbb{P}^n$, embedded to the product in the same way: $\{[z_0: ... : z_n] \times [\bar z_0: ... : \bar z_n] \}$ and then deformed by a 
Hamiltonian isotopy. The arguments are the same as in {\bf Example 1}.

{\bf Example 1b.} Another generalization of {\bf Example 1} is to consider the product $Q_k = \mathbb{C} \mathbb{P}^1 \times ... \times \mathbb{C} \mathbb{P}^1$
with the bundle ${\cal O}(1,...,1)$ where the number of summands is $k> 2$. The computations in this case are much more complicated than in Example 1, 
but several facts can be established in short.

Take $Q_3$ with the bundle ${\cal O}(1,1,1)$ and fix homogenious coordinates $[x_0: x_1], [y_0: y_1], [z_0: z_1]$ on each $\mathbb{C} \mathbb{P}^1$.
Take the section $h_D$ of ${\cal O}(1,1,1)$ given by the polynomial $x_0 y_0 z_0 + x_1 y_1 z_1$. Then the function $\psi^1_D$ has the following finite critical
points and sets: 1) two minimal points at $[1:0] \times [1:0] \times [1:0]$ and (symmetrically) $[0:1]\times [0:1]\times [0:1]$; 
2) critical set $T_c = \{[1: e^{i \eta_x}] \times [1: e^{i \eta_y}] \times [1: e^{i \eta_z}] | \eta_x + \eta_y + \eta_z = 0 \}$ which is isomorphic
to 2 - torus. Taking all gradients lines which start at the minimal points and run to $T_c$ we get a singular submanifolds $S \subset Q_3$,
by the very definition this $S$ is the lagrangian shadow $Sh^{Lag}(D)$ of the irreducible divisor $D = \{ x_0 y_0 z_0 + x_1 y_1 z_1 = 0 \}$. Geometrically this $S$
is isomorphic to 3 - dimensional sphere with two singular points appear after the shrinking of two unknoted smooth circles. Using the toric actions
on the direct summands in $Q_3$ we can extend this observation to generic smooth divisors from the same complete linear system.  On the other hand
for the reducible divisors the lagrangian shadows vanish. It follows from the fact that if $D$ is reducible then $h_D$ is the tensor product of
holomorphic sections of ${\cal O}(1,1,0)$ and ${\cal O}(0,0,1)$ up to transpositions of 1's, therefore the function $\phi^1_D$ is the sum of two functions
and it does admit no finite critical points of Morse index 3. But as it was explained in [1] for the existence of special Bohr - Sommerfeld cycle
one needs critical points of this index.

The situation in the case of $Q_4$ is much reacher: in this case even for reducible divisors from the complete linear system $\vert h_1+ ... + h_4 \vert$
where $h_i$ is the generator for i'th summand, one has nontrivial lagrangian shadows. Indeed, let us take $h_D = h_{D_1} \otimes h_{D_2}$ where 
$D_1 \in \vert h_1 + h_2 \vert, D_2 \in \vert h_3 + h_4 \vert$ and both $D_i$ are irreducible. Then for each $D_i$ we get, according to {\bf Example 1},
the lagrangian shadows lifted from the first and the second pairs of direct summands in $Q_4$ and totally for $D$ we get $Sh^{Lag}(D)$ isomorphic to
the direct product $S^2 \times S^2$ of 2 - spheres. The same is true, of course, for any other division of the set $\{h_i \}$ into pairs, but 
for the reducible divisors of other types (3+1) we again get vanishing lagrangian shadows. Irreducible divisors shadow essentially the same
as in the case of $Q_3$: the difference is in the structure of two singular points in $Sh^{Lag}(D)$.

Our last example is more geometrically interesting:

{\bf Example 2.} Let $X$ be the full flag variety $F^3$ for $\mathbb{C}^3$ with the standard Kahler structure coming from the realization of $F^3$ as the incidence cycle in the direct
product $\mathbb{C} \mathbb{P}^2 \times \mathbb{C} \mathbb{P}^2$, and the line bundle $L = {\cal O}(1,1)|_X = K^{- \frac{1}{2}}_X$.
Then we claim that for irreducible sections  of the bundle the lagrangian shadows of the corresponding divisors is Hamiltonian equivalent to
lagrangian 3 - sphere which is known from the Gelfand - Zeytlin system consideration, see [4] and references therein. This is also true for certain reducible divisors from the
same linear system.  Our arguments are based on the {\bf Example 1a} above and on some functorial property of Special Bohr - Sommerfeld geometry.

{\bf Proposition}. {\it Let $Y \subset X$ is a smooth algebraic subvariety in $X$, and $D_Y \subset Y$ is the intersection $D \cap Y$ for a very ample divisor
$D \subset X$. Suppose that this intersection is transversal. Then if the lagrangian shadow $Sh^{Lag}(D) \subset X$ intersects $Y$ cotransversally so 
if $\rm{dim}_{\mathbb{R}} (Sh^{Lag}(D) \cap Y = \rm{dim}_{\mathbb{C}} Y$ then this intersection equals to $Sh^{Lag}(D_Y)$ in $Y$.}

The proof is direct: take the corresponding line bundle $L_D \to X$, fix the Kahler form $\omega$, then take the corresponding hermitian connection on
$L_D$. Then restricting all the data on $Y$ we get that $D_Y$ is defined by $h_D|_Y$, and the covariantly constant section of $(L_D, a)|_{Sh^{Lag}(D)}$
restricts to a covariantly constant section of    $(L_D, a)|_{Sh^{Lag}(D)\cap Y}$, so the proportionality coefficient function
used to be real positive after any restrictions, which gives us the specialty condition for the intersection $Sh^{Lag} \cap Y$
with respect to the holomorphic section $h_D|_Y$, so it remains only one important condition for $Sh^{Lag} \cap Y$ --- to be
lagrangian in $Y$, what essentially means that the intersection is {\it cotransversal}.

Now let $X = \mathbb{C} \mathbb{P}^2 \times \mathbb{C} \mathbb{P}^2$, algebraic subvariety $Y = F^3$ is given by the equation $\sum x_i y_i = 0$ for fixed
homogenious coordinates on both $\mathbb{C} \mathbb{P}^2$'s, and the divisor $D$ is given by the polynomial $x_0 y_0 + x_1 y_1 - x_2 y_2 = 0$. Then it is not hard 
to see that $Sh^{Lag}(D)$ is defined by the relation $\{ [x_0: x_1: x_2] \times [\bar x_0: \bar x_1: - \bar x_2] \}$ in $X$, being a copy of the projective
plane. But according to [4] the intersection $Sh^{Lag}(D) \cap Y$ is exactly the Gelfand - Zeytlin lagrangian sphere; and therefore it is the lagrangian shadow
of $D_Y = D \cap Y$. Now the point is that $D_Y$ is reducible in $Y = F^3$: it is formed by two del Pezzo surfaces of degree 8, see [2]. However
we can vary the section $h_{D_Y} \to h_{D_ir}$ taking the polynomial $\alpha_0 x_0 y_0 + \alpha_1 x_1 y_1 - \alpha_2 x_2 y_2$, and if the expression
$$
\frac{ \alpha_0 x_0 y_0 + \alpha_1 x_1 y_1 - \alpha_2 x_2 y_2}{x_0 y_0 + x_1 y_1 - x_2 y_2}|_{Sh^{Lag}(D_Y)}
$$
is real positive everywhere then $Sh^{Lag}(D_{ir}) = Sh^{Lag} (D_Y)$, see [1]. This happens if every $\alpha_i \in \mathbb{R}_+$, and this gives us
big family of irreducible divisors with the same lagrangian shadow.

The arguments from the toric geometry lead to the claim that the same is true for generic irreducible divisor from the complete linear system $\vert (h_1 + h_2)|_{F^3} \vert
= \vert - \frac{1}{2} K_{F^3} \vert$.

We end these notes with the following final remark\footnote{I would like to cordially thank Stefan Nemirovsky for the discussion on this
question}. It was conjectured in [1] that under the variation of ample divisors in the complete
 linear system $\vert D \vert$ the ``lagrangian shadows'' keep stable: the ``number'' of components in $Sh^{Lag}(D)$ is the same.
This conjecture is very naive, so  in general setup it is much more reasonable to formulate such a stability property  on the homological level.
Since even in the best case a ``lagrangian shadow''comes without any preffered orientation, we are speaking about
$\rm{mod} \mathbb{Z}_2$ vesrion: the collection of components of a lagrangian shadow gives a class from $H_n(X, \mathbb{Z}_2)$.
It seems that this class must be the same for generic elements of $\vert D \vert$. It implies that in the case when $D$ is very ample (so $k=1$) this class
depends on the Kahler form $\omega_1$ only. Does this class depend on the cohomology class of $\omega_1$ rather than on the form itself? From the Hard Lefschetz
we know that the cohomology class $[\omega_1]$ gives  isomorphisms $H^{n-k} \to H^{n+k}$ but doesn't touch the middle part $H^n$, so may be
there  is some hidden part of the Hard Lefschetz theorem which concerns the middle cohomology group, and the story with lagrangian shadows just reflects
this one.

 $$$$
{\bf References}

[1] Nik. A. Tyurin, ``Special Bohr - Sommerfeld geometry'',  arXiv:1508.06804;

[2] P. Griffits, J. Harris, ``Principles of algrebraic geometry'', NY, Wiley, 1978;

[3] Ya. Eliashberg, ``Topological characterization of Stein manifolds of $\rm{dim} >2$'', Internat. J. Math., 1, no. 1, pp. 29 -46 (1990);

[4] Yu. Nohara, K. Ueda, ``Floer cohomologies of non - toric fibers of the Gelfand - Zeytlin systems'', arXiv:1409.4049.

\end{document}